**Rethinking how randomness is taught in K-12 statistics subjects**

Mark Louie F. Ramos
Department of Health Policy and Administration, Penn State University, University Park, PA

Abstract

The effective teaching and learning of statistics persist as a challenge in K-12 education and has clear impacts in developing competence and confidence of students in entering STEM fields especially in today's digital age of data science. In this paper, it is proposed that the current pedagogy around teaching the concept of randomness, a core element of statistics, is incomplete. It is argued that expanding on the observer-relative nature of randomness will address commonly existing ambiguities on what is and what is not random, help in explaining more complex statistical concepts such as hypothesis testing, and impart greater awareness on the practical value of the practice of statistics in data science.

1. **INTRODUCTION**

It is no mystery that higher order thinking skills are built upon lower-level conceptual and practical knowledge (Maslihah et al. 2020) and that both positive attitudes towards and performance in Science, Technology, Engineering, and Mathematics (STEM) college fields are strongly predicted by students' competence and confidence in mathematics that is developed from elementary through high school (Daker et al. 2021; NCTM 2018; Nitzan-Tamar and Kohen 2022). In today's digital age, data science has emerged as a key area of interest across different STEM fields, emphasizing the critical role of developing statistical literacy alongside other mathematics skills (Watson et al. 2020). Crucial to statistical literacy is teaching the concept of randomness and how it is used to justify practically every statistical tool, yet this is a persistent challenge (Ingram 2022). In a study on student teachers' conceptions of randomness and probability, it was found that even student teachers with a strong background in mathematics may struggle with the concept of randomness and how to teach it (Ingram 2022). In this paper, I discuss the progress of teaching statistics in the K-12 curriculum in the United States and argue that current pedagogy in explaining the concept of randomness is incomplete and needs to be revised in order to dispel misconceptions and provide accurate context in explaining how and why methods of statistical analysis work.

2. **BACKGROUND**

2.1. **Status quo of "random" in K-12 classrooms**

The United States started introducing statistics as a basic high school requirement in the 1920's, with varying levels of progress across the decades that followed (Scheaffer and Jacobbe 2014). In 1989, the National Council of Teachers of Mathematics (NCTM) released curriculum standards that included objectives for K-12 on selecting appropriate statistical methods for analyzing data and understanding basic concepts of statistics and probability (NCTM 1989). By 2014, the Common Core State Standards in Mathematics (CCSS-M) included that by $7^{th}$ grade, students are taught about randomized data and inference (Scheaffer and Jacobbe 2014), and it is notable that the CCSS-M mentioned random in the context of random sampling, e.g., from a finite population, or random assignment, e.g., into treatment and control groups, but not in the context of random phenomenon e.g., what does it mean when something occurs at random or how do we differentiate between something that occurs at random and something that does not? More recently, the Guidelines for Assessment and Instruction in Statistics Education II (GAISE II) was released by the American Statistical Association (ASA) and endorsed by the NCTM to provide updated guidelines to developing students' skills for making sense of data (Bargagliotti et al. 2020). These guidelines include extensive content on randomness and its importance in the practice of data science. Despite of this, issues on the implementation of effective statistics education curricula persist. Batanero (2015) conducted an extensive review of the treatment of randomness as a concept in modern day educational curricula and found that inconsistencies in its definition across curricular documents were common, and that research suggests persistent misconceptions among students. Similarly, Watson and Fitzallen (2019) found that mathematics curriculum documents in various countries including the United States use the word "random" in their learning outcomes but do not provide a definition for the term. (Ingram 2022) claimed that in mathematics education, there is no universally accepted definition of randomness. Furthermore, they detailed how teachers struggled with content knowledge in probability and statistics, and how this has led to an avoidance in examining important concepts such as randomness sufficiently in their pedagogy (Ingram 2022). Khan Academy (Kunz 2021), a non-profit educational content developer that is partnered with over 280 school districts in the U.S., introduces the concept of randomness in Unit 9 of its Probability and Statistics course in the context of random variables. Specifically, it defines random variables as "ways to map outcomes of random processes to numbers," which is accurate but leaves out the definition of a random process and instead substitutes a definition by examples of coin flipping, die rolling, and measuring rainfall (Khan Academy 2012). Even in the comments section of this video, we can see the confusion of viewers over what is actually random; responses from other commenters do not provide any clarity and in some cases, add to the confusion. The rarity of a formal definition seems to stem from ambiguity between the use of the word random in the colloquial sense and its use in mathematics in a technical sense. As pointed out by Watson and Fitzallen (2019), students may think of random as an adjective to describe an unlikely, inexplicable, or surprising event.

### 2.2. So what is randomness?

Randomness is formally and accurately defined in various sources and in different phrasing as uncertainty in the outcomes of experiments, activities, or observations (Franklin and Bargagliotti 2020; Pishro-Nik 2014; Watson and Fitzallen 2019). Moreover, Franklin and Bargagliotti (2020) go on to discuss how key to understanding statistical concepts is to develop probabilistic thinking, which involves thinking about variability in the data that we observe. Batanero(2015) provides an exhaustive history of the definition of randomness from primitive ideas to modern day educational curricula. A common aspect of these definitions that is critical but is hidden in the background is the fact that randomness is relative. Spiegelhalter (2024) wrote about how probability, the way we measure randomness as a number between 0 and 1, "does not exist." This was to say that probability is not some real quantity, like temperature or force, but is a construct with values that are entirely dependent on the knowledge of the perceiver. For example, the boiling point of water exists regardless of there being anyone to perceive it. The boiling point is not constant since it changes, for example, with altitude, but it exists independent of perception. In contrast, the probability that a specific kettle of water will start boiling when it reaches 100°C is not independent of perception. Whether or not it will start boiling is already determined by the physical factors surrounding the activity. Thus, in this context, the probability is simply 1 if it will happen and 0 if it will not. Now, consider if there are people observing this event and need to estimate the probability. Those who have absolutely no knowledge about the topic of boiling points may think the probability is 0.50, while those who learned that 100°C is the boiling point of water but fail to account for altitude may think that the probability is 1. Those who know that boiling points are relative to altitude may not give a probability of 1 because they are considering how far below sea level the kettle is at, but they may still say the probability is close to 1. Meanwhile, the person who knows that the kettle is in Death Valley, California which is 282 feet below sea level may give a probability of 0, because they know that this is a low enough altitude to increase the boiling point. While these examples are entirely from a Bayesian perspective of computing probability, we can also add the frequentist perspective. A person in this example who is a frequentist statistician would want to try and see how many of 100 identical kettles of water would boil at 100°C. Suppose another person who has already done this experiment before tells him that only 10 out of 100 did so, then for the frequentist, the estimated probability is 0.10. None of these are the "right" probability; none of them are the "wrong" probability either, because they are all measurements of randomness that are relative to what the observer knows and what rules they decide to use.

Simply put, something is random because we, the observers, lack the ability to predict its outcomes perfectly. If we could, then it is no longer random to us, but may still be random to those who could not. Franklin and Bargagliotti(2020) allude to this when they discussed how the models in statistics are based on assumptions that are made in the context of the situation. They say that these assumptions may not stay constant, implying that they may change once the data is seen. This reflects the relative nature of randomness; statisticians make the best assumptions that they can make at the start of the research given how much of the data generating process they understand and how much they do not.

Various works in statistics education emphasize helping students understand the concept of randomness. GAISE II recommends for students to experience randomness in dynamic settings via experimentation or simulation (Bargagliotti et al. 2020), and this has been incorporated as part of learning expectations in many K-8 schools in the United States (Weiland and and Sundrani 2022). Various recommended strategies, such as those by Watson and Fitzallen (2019) or Batista et al. (2022) focus on developing student intuition about uncertainty and variation through experience. While these strategies are useful, they do not directly address the observer-relative nature of randomness that is focused on in this paper. For example, the GAISE II recommends for students to compare theoretical probability with experimental probability, such as by having students toss a coin multiple times or be shown a simulation of such to see that the experimental probability converges to the theoretical probability with enough trials. The observer-relative idea of randomness in this paper asks students and teachers to take a step back before this and understand that this exercise is only useful because one cannot tell how the coin lands each time. If one could, then the concept of probability becomes practically useless for that person. If you know that a coin will land heads with certainty, then for you the probability that it will do so is 1.

## 3. PROPOSED LEARNING CONTENT CHANGES WHEN EXPLAINING RANDOMNESS

### 3.1. Conceptual Grounding

This section is prefaced with the caveat that the proposed learning content changes are based on the idea that randomness is relative to the observer; they have not yet been formally tested in classroom experiments to measure their effectiveness. Instead, they are grounded on more than a decade of teaching experience in introductory statistics courses from high school to undergraduate levels and, perhaps more importantly, repeated exposure to these courses as a student, from high school to graduate school. As a student, learning statistics in high school, all the way through an undergraduate mathematics program and, later, a master's in mathematics education program, mainly felt like a repetition on determining which test went with which situation, an experience that is similarly reflected in various works (Gigerenzer 2004; Lee 1997). It was not until experiencing statistics taught with a focus on understanding its probabilistic underpinnings, rather than these being just briefly mentioned then forgotten, that I could make a sincere appreciation of its methods and principles. This led to changes in how I taught introductory statistics courses, with some specific examples as described in Ramos (2025). Making these changes helped me better convey to my students both how and why statistical tools work. Like I did when I first understood these as a student, my own students appreciated being explicitly told that randomness is not some constant quantity but is individually relative to what they know or assume. These changes are simple and straightforward to implement and are therefore offered as an invitation for readers to consider incorporating them into their pedagogy if they find the rationale compelling.

### 3.2. The Omniscient Observer

The foremost change is emphasizing how randomness is not just tied to uncertainty, but that uncertainty is subjective. To do this, teachers can introduce the concept of the **omniscient observer** at the beginning of the lesson on defining probability. The omniscient observer is someone who has all the information necessary to know what will happen before it happens. Such an observer does exist in specific contexts. For example, for coin tossing, these two videos show examples of people who built machines that toss a coin the same way each time (Briggs 2023; Lenzonhighst 2015), and thus they act as omniscient observers who know how the coin will land prior to tossing it. This demonstrates that if someone has all the tools needed for an observation, then they have perfect information on the outcome of the event. From here, the teacher can pivot to what the typical coin toss looks like and explain that suppose humans could measure every physical variable that goes into a coin toss, then we would be omniscient observers of this scenario. However, since that is not the case, then we create the concept of randomness to reflect our uncertainty about the things we observe. This serves as a good springboard for introducing the definition of probability, a quantity with which we measure or express the randomness in a situation. In the example of tossing a coin, teachers would typically introduce the idea of equiprobability; a coin has two outcomes, heads and tails, and so the probability of heads is ½. Yet this is only true if each outcome is equally likely, of which we are not certain. Therefore, we only assume that they are equally likely and use that as basis for our probability. Applying this to 6-sided die rolls, without any other information, we can assume that each outcome is equally likely and therefore say that the probability of each one is 1/6. This can in turn be used as a springboard for introducing the frequentist definition of probability, where you can base it from observing the event of interest a large number of times. The teacher can show a simulation or a frequency table that shows how some values in the die appear considerably more times than the others over many throws. This lets us refine out measurement of probability for this particular example beyond our initial assumption of equiprobability and conclude that the die is loaded. The key change from current pedagogy is establishing the idea that these "rules" for how to compute probabilities are grounded on assumptions that are subjective to our uncertainty (Spiegelhalter 2024). That is, these rules allow us to formalize the entire body of mathematics governing the computation of probabilities, but the rules themselves are grounded on subjective assumptions. For example, if we assume that a group of 10 coin tosses are independent of each other and each coin toss is fair, then that allows us to compute the probabilities for the number of heads in 10 tosses based on a binomial distribution with parameters $n=10$ and $p=0.50$.

### 3.3. The Infallible Statistician

As described, the main change that is proposed in this paper is simple and easily implementable. In terms of time investment, it amounts to no more than a fraction of a class period – 15 minutes of a single, 60-minute class session and perhaps 3-minute reminders in succeeding sessions when relevant. However, making this change sets a tone for the rest of the course that the teacher must be mindful of when introducing more complex concepts later. For example, it should be pointed out that the definitions of random variable and random distribution do not change in the presence of an omniscient observer, it is just that such an observer has no practical use for those concepts. The random variable associated with a coin toss is typically defined as X with a support of {0,1}, where 0 and 1 are arbitrarily assigned to tails and heads respectively. If the coin tossing process is fair, then $Pr(X=1)=Pr(X=0)=0.5$ which completely defines the random variable X and its distribution. The omniscient observer may not necessarily have any control over how the coin lands, so if we observe this process repeatedly with the omniscient observer, both of us will observe the same phenomenon of the coin landing heads about half the time. The difference is that the omniscient observer can tell how the coin will land each time, while we cannot. As such, this reinforces why we need to use the definitions of random variable and random distribution to help us make decisions related to the experiment. For example, suppose that betting on heads pays $1 while betting on tails pays $1.10 at the same cost. If the coin is fair then the best strategy for us is to always bet tails. In contrast, the omniscient observer has no use for that strategy, since they can choose the winning outcome each time. This brings up another important concept, which is that of the **infallible statistician**. The infallible statistician is one whose assumptions are always correct. We note that our strategy to bet on heads each time is based on our assumption that the coin is fair. If we were infallible statisticians, then that strategy would be the best strategy. However, the truth is that just like we are practically never omniscient observers, we are also practically never infallible statisticians. If the truth was that the probability of heads was actually 0.53 instead of 0.50, then the optimal strategy would shift to always betting heads at the same payoff values. As fallible statisticians, we *choose* to assume that the coin is fair unless and until we see evidence to the contrary, and then we can choose to change our assumption, but no practical amount of evidence can provide us with absolute certainty that even that change is correct. We argue that being both cognizant and comfortable with this reality of subjective uncertainty is important to teaching and learning statistics.

## 4. EDUCATIONAL RELEVANCE OF PROPOSED CHANGES

By adopting the proposed changes, we remove the ambiguity about what is or is not random in discussions that happen in the K-12 classroom. The confusion that arises in defining what is random is caused by missing the idea that with enough information, nothing is actually random. If students understood this from the beginning, then they will know, for example, that assigning an equiprobable distribution to a coin toss or a die roll is a choice, not a law. We do make this choice in our lessons which is instructive so that the probability rules we teach abide by these assumptions and we all get the same answers for specific problems, but the practical reality does not necessarily abide by any of the assumptions we make in the classroom.

Moreover, by bringing this up early and reinforcing it as needed in teaching related concepts, we help students build intuition that is key to understanding how and why statistical tools work.

George Box famously said, "All models are wrong, some are useful" (Box 1976). This underscores how all model selections made in the practice of statistics are approximations based on observation, analytical reasoning, substantive knowledge, or all three rather than the unerring declarations that they can seem to be, especially to impressionable high school students being exposed to these concepts for the first time. For example, in simple linear regression, the basic model is introduced as $y = \beta_0 + \beta_1 x + \varepsilon$, where $\beta_0$ and $\beta_1$ are fixed but unknown parameters and randomness is introduced only through the error term $\varepsilon$. This model choice for any practical application is almost certainly wrong. The true relationship between $y$ and $x$ could in reality be $y = \beta_0 + \beta_1 x^{1.1} + \varepsilon$ or something entirely different. The idea that we, being fallible statisticians, can pick the exact right one among an infinity of possible choices is bluntly absurd. Yet despite of this, the practical utility of regression models is widely evident. That is, we do not need to get the model assumption exactly correct for the model to be useful. For example, the relationship between how much a person smokes and whether or not they develop lung cancer was established based on similarly simple logistic regression models in landmark epidemiological studies during the 1930's to the 1950's (Thun 2010). These studies, despite using flawed models, were responsible for critical public health initiatives credited for reducing the disease's prevalence decades later (Thun 2010). Taking this a step further, core to the idea of statistical inference (statistical hypothesis testing) is making a decision after observing the test statistic produced under the assumption that the null distribution is true. Rejection of the null hypothesis is in essence an admission that the initial assumption made about the distribution where the data came from is dubious enough to reasonably justify acting as though it is false. Specifying an actual probability for the null hypothesis being false is not essential to making this decision. This only makes sense in the underlying notion of randomness being relative to the observer and helps to explain why the argument that p-values do not provide a probability about the null hypothesis, while absolutely true, is at best a strawman criticism of the p-value's utility when used properly.

The statistician makes assumptions and selects distributions based on careful study and professional experience, often in collaboration with multidisciplinary research partners, but they **must** do so because they are neither omniscient observers nor infallible statisticians. That is, their lack of knowledge is consciously modeled as randomness in the data and applying statistics gets to answer questions about the process that produced the data within some degree of confidence. Furthermore, the proposed pedagogical change imparts lessons that can make it easier to communicate issues in the practice of statistics that arise in higher level data science courses and even more so, in actual data analysis practice. Many practical issues in data collection, such as nonrandom missing data and participant response rates, as well as issues on measurement error or generalization of results to other populations, are tied to the notion that randomness is relative to the observer, and that the researcher's inability to know different aspects of the data generating process is what creates these problems that need to be addressed. For example, in dealing with missing data, assuming that the data are "missing completely at random" (MCAR), makes the optimal solution for data imputation straightforward (van Buuren 2018), yet such a strong assumption is often difficult if not expressly incorrect to make. Inculcating the relative nature of randomness can help prepare the students' mindset to understand why data imputation is a complicated problem and how incorrect solutions can lead to faulty scientific conclusions (Alwateer et al. 2024).

## 5. CONCLUSIONS AND FUTURE DIRECTION

We need to change how we teach the concept of randomness in high school statistics. By putting some elaboration on the importance of "we" when defining random phenomena as processes that yield outcomes which we, as humans, have a level of uncertainty over, misconceptions on what is and what is not random are dispelled. Students will understand that randomness and probability are not natural constructs but are purely manmade concepts that depend entirely on how much about the phenomenon we are uncertain about. By doing this, we fill a knowledge gap and build intuition that can help students understand more complex statistical concepts, particularly when it comes to how and why statistical tools work in practice.

Beyond anecdotal evidence from personal experience, the ideas presented in this paper have not yet been formally evaluated. A natural next step is to design quasi-experimental studies or, where possible, randomized controlled experiments, that test their effectiveness in real classrooms. While we plan to pursue such studies, a persistent challenge in education research is the large heterogeneity of effects across settings. This makes it valuable for individual statistics instructors who find these ideas worthwhile to conduct their own pilot studies of these strategies in their own classrooms.

**Acknowledgements**

I would like to thank my mentor, Barry Graubard, Ph.D., Scientist Emeritus at the National Cancer Institute, for his valuable insights into this paper. I also thank the editors and reviewers of the Journal of Statistics and Data Science Education for their effort and expertise in reviewing my manuscript. For transparency, this paper was submitted in December 28, 2024 and underwent 3 rounds of reviews with mainly positive feedback. I withdrew it last February 15, 2026 because I did not have the resources to conduct empirical studies which the editor had suggested and considered necessary for potential publication.

**Funding Statement**

There is no funding to declare.

**Disclosure**

There are no conflicts of interest to declare.

**Data Availability Statement**

Data sharing is not applicable. There were no primary or secondary datasets collected or analyzed.

**Generative AI Statement**

Generative AI was not used in writing this paper.

# References

Alwateer, M., Atlam, E.-S., El-Raouf, M. M. A., Ghoneim, O. A., and Gad, I. (2024), "Missing Data Imputation: A Comprehensive Review," *Journal of Computer and Communications*, Scientific Research Publishing, 12, 53–75. https://doi.org/10.4236/jcc.2024.1211004.

Bargagliotti, A., Franklin, C., Arnold, P., Gould, R., Johnson, S., Perez, L., and Spangler, D. (2020), "Guidelines for Assessment and Instruction in Statistics Education (GAISE) Reports," Available at https://www.amstat.org/education/guidelines-for-assessment-and-instruction-in-statistics-education-(gaise)-reports.

Batanero, C. (2015), "Understanding Randomness: Challenges for Research and Teaching," Prague, Czech Republic.

Batista, R., Borba, R., and Henriques, A. (2022), "Fairness In Games: A Study on Children's and Adults' Understanding of Probability," *Statistics Education Research Journal*, 21, 13–13. https://doi.org/10.52041/serj.v21i1.79.

Box, G. E. P. (1976), "Science and Statistics," *Journal of the American Statistical Association*, ASA Website, 71, 791–799. https://doi.org/10.1080/01621459.1976.10480949.

Briggs, W. (2023), "The probability of a head in a coin flip is 1, not 1/2!," Available at https://www.youtube.com/watch?v=8XX0iRAN--8.

van Buuren, S. (2018), *Flexible Imputation of Missing Data*, CRC Press.

Daker, R. J., Gattas, S. U., Sokolowski, H. M., Green, A. E., and Lyons, I. M. (2021), "First-year students' math anxiety predicts STEM avoidance and underperformance throughout university, independently of math ability," *NPJ Science of Learning*, 6, 17. https://doi.org/10.1038/s41539-021-00095-7.

Franklin, C., and Bargagliotti, A. (2020), "Introducing GAISE II: A Guideline for Precollege Statistics and Data Science Education," *Harvard Data Science Review*, The MIT Press, 2. https://doi.org/10.1162/99608f92.246107bb.

Gigerenzer, G. (2004), "Mindless statistics," *The Journal of Socio-Economics*, Statistical Significance, 33, 587–606. https://doi.org/10.1016/j.socec.2004.09.033.

Ingram, J. (2022), "Randomness and probability: exploring student teachers' conceptions," *Mathematical Thinking and Learning*, Routledge, 26, 1–19. https://doi.org/10.1080/10986065.2021.2016029.

Khan Academy (2012), "Random variables," Available at https://www.khanacademy.org/math/statistics-probability/random-variables-stats-library/random-variables-discrete/v/random-variables.

Kunz, B. (2021), "More than 280 School Districts Enroll in Khan Academy, NWEA Mastery-Based Learning Offerings," *Khan Academy Blog*, Available at https://blog.khanacademy.org/more-than-280-school-districts-enroll-in-khan-academy-nwea-mastery-based-learning-offerings/.

Lee, P. M. (1997), *Bayesian Statistics: An Introduction*, London: Hodder Education Publishers.

Lenzonhighst (2015), "Bill V's coin flipping machine," Available at https://www.youtube.com/watch?v=lGO72E1V-3E.

Maslihah, S., Waluya, S. B., Rochmad, and Suyitno, A. (2020), "The Role Of Mathematical Literacy To Improve High Order Thinking Skills," *Journal of Physics: Conference Series*, IOP Publishing, 1539, 012085. https://doi.org/10.1088/1742-6596/1539/1/012085.

NCTM (1989), *Curriculum and evaluation standards for school mathematics*, Reston, Va. : The Council.

NCTM (2018), "Building STEM Education on a Sound Mathematical Foundation - National Council of Teachers of Mathematics," Available at https://www.nctm.org/Standards-and-

Positions/Position-Statements/Building-STEM-Education-on-a-Sound-Mathematical-Foundation/.

Nitzan-Tamar, O., and Kohen, Z. (2022), "Secondary school mathematics and entrance into the STEM professions: a longitudinal study," *International Journal of STEM Education*, 9, 63. https://doi.org/10.1186/s40594-022-00381-9.

Pishro-Nik, H. (2014), *Introduction to Probability, Statistics, and Random Processes*, Kappa Research, LLC.

Ramos, M. L. (2025), "On the Pedagogy of Randomness," *Statistics Teacher*.

Scheaffer, R. L., and Jacobbe, T. (2014), "Statistics Education in the K-12 Schools of the United States: A Brief History," *Journal of Statistics Education*, 22, 8. https://doi.org/10.1080/10691898.2014.11889705.

Spiegelhalter, D. (2024), "Why probability probably doesn't exist (but it is useful to act like it does)," *Nature*, 636, 560–563. https://doi.org/10.1038/d41586-024-04096-5.

Thun, M. J. (2010), "Early landmark studies of smoking and lung cancer," *The Lancet Oncology*, Elsevier, 11, 1200. https://doi.org/10.1016/S1470-2045(09)70401-2.

Watson, J., and Fitzallen, N. (2019), "Building Understanding of Randomness from Ideas about Variation and Expectation," *Statistics Teacher*.

Watson, J., Fitzallen, N., and Chick, H. (2020), "What Is the Role of Statistics in Integrating STEM Education?," in *Integrated Approaches to STEM Education: An International Perspective*, eds. J. Anderson and Y. Li, Cham: Springer International Publishing, pp. 91–115. https://doi.org/10.1007/978-3-030-52229-2_6.

Weiland, T., and and Sundrani, A. (2022), "Opportunities for K-8 Students to Learn Statistics Created by States' Standards in the United States," *Journal of Statistics and Data Science Education*, Taylor & Francis, 30, 165–178. https://doi.org/10.1080/26939169.2022.2075814.